\documentclass[a4paper,11pt]{amsart}
\usepackage{graphicx}
\usepackage{mathptmx}
\usepackage{amsmath}
\usepackage{amssymb}
\usepackage{enumitem}
\usepackage{xcolor}

\newmuskip\pFqmuskip

\newcommand*\pFq[6][8]{%
  \begingroup 
  \pFqmuskip=#1mu\relax
  \mathcode`=\string"8000
  \begingroup\lccode`\~=`\,
  \lowercase{\endgroup\let~}\pFqcomma
  F^{#2}_{#3}{\left(\genfrac..{0pt}{}{#4}{#5}\bigg|#6\right)}%
  \endgroup
}
\newcommand{\pFqcomma}{\mskip\pFqmuskip}

\newtheorem{theorem}{Theorem}

\begin{document}

\title[]{Multi-Stirling numbers of the second kind}

\author{Taekyun  Kim$^{1,*}$}
\address{$^{1}$Department of Mathematics, Kwangwoon University, Seoul 139-701, Republic of Korea}
\email{tkkim@kw.ac.kr}

\author{Dae San Kim$^{2,*}$}
\address{$^{2}$Department of Mathematics, Sogang University, Seoul 121-742, Republic of Korea}
\email{dskim@sogang.ac.kr}

\author{Hye Kyung  Kim$^{3,*}$}
\address{$^{3}$Department of Mathematics Education, Daegu Catholic University, Gyeongsan 38430, Korea}
\email{hkkim@cu.ac.kr}

\thanks{* is corresponding author}

\begin{abstract}
The multi-Stirling numbers of the second kind, the unsigned multi-Stirling numbers of the first kind, the multi-Lah numbers and the multi-Bernoulli numbers are all defined with the help of the multiple logarithm, and generalize respectively the Stirling numbers of the second kind, the unsigned Stirling numbers of the first kind, the unsigned Lah numbers and the higher-order Bernoulli numbers .
The aim of this paper is to introduce the multi-Stirling numbers of the second kind and to find several identities involving those four numbers defined by means of the multiple logarithm and some other special numbers.

\end{abstract}

 \maketitle

\markboth{\centerline{\scriptsize  Multi-Stirling numbers of the second kind}}
{\centerline{\scriptsize   T. Kim, D. S. Kim, and H. K. Kim}}

\section{Introduction}

\medskip

The Stirling number of the second kind                                                                                                                                                                                                                                                       ${n\brace r}$ counts the number of ways to partition a set with $n$ elements into $r$ non-empty subsets. These numbers are generalized to the multi-Stirling numbers of the second kind ${n \brace k_1,k_2,\cdots,k_r}$ (see \eqref{eq18}) which reduce to the Stirling numbers of the second kind for $(k_{1},k_{2},\dots,k_{r})=(1,1,\dots,1)$. Indeed, ${n \brace 1,1,\cdots,1}={n\brace r}$.\par
On the other hand, the unsigned Stirling number of the first kind ${n\brack r}$ enumerates the number of permutations of a set with $n$ elements which are products of $r$ disjoint cycles. These numbers are generalized to the unsigned multi-Stirling numbers of the first kind ${n \brack k_{1},k_{2},\cdots,k_{r}}$ (see \eqref{eq28}) which become the unsigned Stirling numbers of the first kind for $(k_{1},k_{2},\dots,k_{r})=(1,1,\dots,1)$. Thus, ${n \brack 1,1,\cdots,1}={n\brack r}$.\par
The unsigned Lah numbers $L(n,k)$ counts the number of ways of a set of $n$ elements can be partitioned into $k$ nonempty linearly ordered subsets. These numbers are generalized to the multi-Lah numbers $L^{(k_{1},k_{2},\dots,k_{r})}(n,r)$ (see \eqref{eq08}) which boil down to the unsigned Lah numbers for $(k_{1},k_{2},\dots,k_{r})=(1,1,\dots,1)$. Indeed, $L^{(1,1,\dots,1)}(n,r)=L(n,r).$ \par
In addition, the multi-Bernoulli numbers $B_{n}^{(k_1,k_2,\dots,k_{r})}$ (see \eqref{eq07}) generalize the Bernoulli numbers of order $r$. In fact, we see that $B_{n}^{(1,1,\dots,1)}=\frac{1}{r!}(-1)^{n}B_{n}^{(r)}$. \par
The common feature of those four numbers is that they are all defined with the help of the multiple logarithm $\mathrm{Li}_{k_{1},k_{2},\dots,k_{r}}(z)$ (see \eqref{eq04}), which reduce to the poly-logarithm $\mathrm{Li}_{k_{1}}(z)$, for $r=1$. \par
The aim of this paper is to introduce the multi-Stirling numbers of the second kind and to derive several identities involving those four numbers defined by means of the multiple logarithm and some other numbers.
In more detail, in Theorem 3 we derive an identity involving the multi-Bernoulli numbers, the Stirling numbers of the second kind, Bell polynomials and the multi-Stirling numbers of the second kind. In Theorem 4, we express the multi-Stirling number of the second kind in terms of the unsigned multi-Stirling numbers of the first kind. We find in Theorem 5 an identity involving the Stirling numbers of the second kind, the multi-Lah numbers, the multi-Stirling numbers of the second kind and the Fubini polynomials of order $r$. We deduce in Theorem 6 an identity involving the multi-Stirling numbers of the second kind, the signed Stirling numbers of the first kind and the multi-Lah numbers. Finally, in Theorem 7 we show an identity involving the multi-Stirling numbers of the second kind, the signed Stirling numbers of the first kind, the Stirling numbers of the second kind and the multi-Bernoulli numbers.
For the rest of this section, we recall the necessary facts that will be needed throughout this paper.

\vspace{0.1in}
We recall that the falling factorial sequence is given by
\begin{equation}\label{eq01}
\begin{split}
(x)_0=1,\quad (x)_n=x(x-1)\cdots(x-n+1),\ \ (n\geq1),\quad (\text{see [1-5,7-16,18,19]}).
\end{split}
\end{equation}
It is well known that the Stirling numbers of the second kind are defined by
\begin{equation}\label{eq02}
\begin{split}
x^n=\sum_{k=0}^n{n \brace k} (x)_k,\quad (n\geq0), \quad (\text{see \cite{5}}).
\end{split}
\end{equation}
From \eqref{eq02}, we note that the generating function of the Stirling numbers of the second kind is given by
\begin{equation}\label{eq03}
\begin{split}
\frac{1}{k!}(e^t-1)^k=\sum_{n=k}^\infty {n \brace k} \frac{t^n}{n!}, \quad (k\geq0),\quad \ (\text{see \cite{16}}).
\end{split}
\end{equation}

For integers $k_i \ (1\leq i\leq r)$, the multiple logarithm (of index $(k_1,\ k_2,\cdots, k_r)$) is given by
\begin{equation}\label{eq04}
\begin{split}
\mathrm{Li}_{k_1,k_2,\cdots,k_r}(t)=\sum_{0<m_1<\cdots<m_r}\frac{t^{m_{r}}}{m_1^{k_1}m_2^{k_2}\cdots m_r^{k_r}},\quad (\text{see \cite{8,10}}).
\end{split}
\end{equation}
If $r=1, \displaystyle \mathrm{Li}_{k_1}(t)=\sum_{m=1}^\infty\frac{t^m}{m^{k_1}}$ is the polylogarithm.
By \eqref{eq04}, we easily get
\begin{equation}\label{eq05}
\begin{split}
\frac{d}{dt}\mathrm{Li}_{k_1,\cdots,k_r}(t)=\frac{1}{t}\mathrm{Li}_{k_1,k_2,\cdots,k_r-1}(t),\quad (\text{see \cite{9,10,15}}).
\end{split}
\end{equation}
From \eqref{eq05}, we have
\begin{equation}\label{eq06}
\begin{split}
\frac{d}{dt}\mathrm{Li}_{k_1,\cdots,k_{r-1},1}(t)&=\frac{1}{t}\sum_{0<m_1<\cdots<m_{r-1}}\frac{1}{m_1^{k_1}\cdots m_{r-1}^{k_{r-1}}}\sum_{m_r=m_{r-1}+1}^\infty t^{m_r}\\
&=\frac{1}{1-t}\mathrm{Li}_{k_1,k_2,\cdots,k_{r-1}}(t).
\end{split}
\end{equation}

In 2000, Kim-Kim introduced the multi-Bernoulli numbers of index $(k_1,k_2,\dots,k_r)$ which are defined by
\begin{equation}\label{eq07}
\begin{split}
\frac{1}{(1-e^{-t})^r}\mathrm{Li}_{k_1,\cdots,k_r}(1-e^{-t})=\sum_{n=0}^\infty B_n^{(k_1,k_2,\cdots,k_r)}\frac{t^n}{n!},\quad (\text{see \cite{9}}).
\end{split}
\end{equation}
Note that $B_{n}^{(\overbrace{1,1,\cdots,1}^{r - \text{times}})}=\frac{(-1)^n}{r!}B_n^{(r)}=\frac{B_n^{(r)}(r)}{r!}$, where $B_n^{(r)}(x)$ are the Bernoulli polynomials of order $r$ given by
\begin{equation*}
\begin{split}
\bigg(\frac{t}{e^t-1}\bigg)^re^{xt}=\sum_{n=0}^\infty B_n^{(r)}(x)\frac{t^n}{n!},\quad (\text{see [9-16]}),
\end{split}
\end{equation*}
and $B_{n}^{(r)}=B_{n}^{(r)}(0)$ are the Bernoulli numbers of order $r$. \par
Recently, Kim-Kim-Kim-Lee-Park introduced the multi-Lah numbers defined by
\begin{equation}\label{eq08}
\begin{split}
\frac{\mathrm{Li}_{k_1,\cdots,k_r}(1-e^{-t})}{(1-t)^r}=\sum_{n=r}^\infty L^{(k_1,k_2,\cdots,k_r)}(n,r)\frac{t^n}{n!},\quad (\text{see \cite{8}}).
\end{split}
\end{equation}
Note that
\begin{equation}\label{eq09}
\begin{split}
L^{(\overbrace{1,1,\cdots,1}^{r- \text{times}})}(n,r)=L(n,r),\quad (n,r\geq0),
\end{split}
\end{equation}
where $L(n,r)$ are the unsigned Lah numbers given by
\begin{equation}\label{eq10}
\begin{split}
\frac{1}{r!}\bigg(\frac{t}{1-t}\bigg)^r=\sum_{n=r}^\infty L(n,r)\frac{t^n}{n!},\quad (\text{see \cite{5,8,16}}).
\end{split}
\end{equation}

Y. Ma et al also studied the multi-Stirling numbers of the first kind defined by
\begin{equation}\label{eq11}
\begin{split}
(-1)^r\mathrm{Li}_{k_1,k_2,\cdots,k_r}(-t)=\sum_{n=r}^\infty S^{(k_1,k_2,\cdots,k_r)}(n,r)\frac{t^n}{n!},\quad (\text{see \cite{15}}).
\end{split}
\end{equation}
Note that $S^{(\overbrace{1,1,\cdots,1}^{r-\text{times}})}=S_1(n,r)$, where $S_1(n,r)$ are the Stirling numbers of the first kind given by
\begin{equation}\label{eq12}
\begin{split}
(x)_n=\sum_{k=0}^nS_1(n,k)x^k,\quad (n\geq0), \quad (\text{see \cite{5, 12}}).
\end{split}
\end{equation}

The unsigned Stirling numbers of the first kind are defined by
\begin{equation}\label{eq13}
\begin{split}
\langle x\rangle_n=\sum_{k=0}^n{n \brack k}x^k,\quad (n\geq0), \quad (\text{see \cite{5, 11}}),
\end{split}
\end{equation}
where $\langle x\rangle_0=1,\ \langle x\rangle_n=x(x+1)\cdots(x+n-1), \ (n\geq1)$.
From \eqref{eq12} and \eqref{eq13}, we note that
\begin{equation}\label{eq14}
\begin{split}
{n \brack k }=(-1)^{n-k}S_1(n,k),\quad (n,k\geq0),\quad (\text{see \cite{5, 11}}).
\end{split}
\end{equation}

Now, we recall that the Fubini polynomials of order $r$ are given by
\begin{equation}\label{eq15}
\begin{split}
\bigg(\frac{1}{1-x(e^t-1)}\bigg)^r=\sum_{n=0}^\infty F_n^{(r)}(x)\frac{t^n}{n!},\quad (\text{see \cite{12}}).
\end{split}
\end{equation}
From \eqref{eq15}, we note that
\begin{equation}\label{eq16}
\begin{split}
\sum_{n=0}^\infty F_n^{(r)}(x)\frac{t^n}{n!}&=\bigg(\frac{1}{1-x(e^t-1)}\bigg)^r=\sum_{k=0}^\infty\binom{r+k-1}{k}x^k(e^t-1)^k\\
&=\sum_{k=0}^\infty x^k\binom{r+k-1}{k}k!\sum_{n=k}^\infty {n \brace k}\frac{t^n}{n!}\\
&=\sum_{n=0}^\infty\bigg(\sum_{k=0}^n \binom{r+k-1}{k}{n \brace k}k!x^k\bigg)\frac{t^n}{n!}.
\end{split}
\end{equation}
Thus, by \eqref{eq16}, we get
\begin{equation}\label{eq17}
\begin{split}
F_n^{(r)}(x)=\sum_{k=0}^n{n \brace k} k!\binom{r+k-1}{k}x^k, \quad (\text{see \cite{12}}).
\end{split}
\end{equation}

\section{Multi-Stirling numbers of the second kind}

In view of \eqref{eq03}, we may consider the multi-Stirling numbers of the second kind defined by
\begin{equation}\label{eq18}
\begin{split}
\mathrm{Li}_{k_1,k_2,\cdots,k_r}(1-e^{1-e^t})=\sum_{n=r}^\infty {n \brace k_1,k_2,\cdots,k_r}\frac{t^n}{n!}.
\end{split}
\end{equation}

Now, we observe that
\begin{equation}\label{eq19}
\begin{split}
\frac{d}{dt}&\mathrm{Li}_{k_1,k_2,\cdots,k_{r-1},1}(1-e^{1-e^t})=\frac{d}{dt}\sum_{0<m_1<\cdots<m_r}\frac{(1-e^{1-e^t})^{m_r}}{m_1^{k_1}\cdots m_{r-1}^{k_{r-1}}m_r}\\
&=\sum_{0<m_1<\cdots<m_{r-1}}\frac{1}{m_1^{k_1}\cdots m_{r-1}^{k_{r-1}}}\frac{e^t e^{1-e^t}}{1-e^{1-e^t}}\sum_{m_r=m_{r-1}+1}^\infty (1-e^{1-e^t})^{m_r}\\
&=e^t\sum_{0<m_1<\cdots<m_{r-1}}\frac{(1-e^{1-e^t})^{m_{r-1}}}{m_1^{k_1}m_2^{k_2}\cdots m_{r-1}^{k_{r-1}}}=e^t\mathrm{Li}_{k_1,\cdots,k_{r-1}}(1-e^{1-e^t}).
\end{split}
\end{equation}

Thus, by \eqref{eq19}, we get
\begin{equation}\label{eq20}
\begin{split}
\frac{d}{dt}\mathrm{Li}_{1,1}(1-e^{1-e^t})=e^t\mathrm{Li}_1(1-e^{1-e^t})=-e^t\log(e^{1-e^t})=e^t(e^t-1).
\end{split}
\end{equation}

From \eqref{eq20}, we note that
\begin{equation}\label{eq21}
\begin{split}
\mathrm{Li}_{1,1}(1-e^{1-e^t})=\int_0^t\mathrm{Li}_{1,1}(1-e^{1-e^x})dx=\int_0^te^x(e^x-1)dx=\frac{1}{2!}(e^t-1)^2.
\end{split}
\end{equation}

By \eqref{eq19} and \eqref{eq20}, we get
\begin{equation}\label{eq22}
\begin{split}
\mathrm{Li}_{1,1,1}(1-e^{1-e^t})&=\int_0^t e^x\mathrm{Li}_{1,1}(1-e^{1-e^x})dx\\
&=\frac{1}{2!}\int_0^te^x(e^x-1)^2dx=\frac{1}{3!}(e^t-1)^3.
\end{split}
\end{equation}

Continuing this process, we have
\begin{equation}\label{eq23}
\begin{split}
\mathrm{Li}_{\underbrace{1,1,\cdots,1}_{r-\text{times}}}(1-e^{1-e^t})=\frac{1}{r!}(e^t-1)^r=\sum_{n=r}^\infty {n \brace r}\frac{t^n}{n!}.
\end{split}
\end{equation}

Therefore, by \eqref{eq18} and \eqref{eq23}, we obtain the following theorem.

\begin{theorem}

For $n,r\geq0$ with $n\geq r$, we have
\begin{equation*}
\begin{split}
{n \brace \underbrace{1,1,\cdots,1}_{r-\text{times}}}={n \brace r}.
\end{split}
\end{equation*}
\end{theorem}

From \eqref{eq19}, we note that
\begin{equation}\label{eq24}
\begin{split}
\frac{d}{dt}&\mathrm{Li}_{k_1,k_2,\cdots,k_{r-1},1}(1-e^{1-e^t})=e^t\mathrm{Li}_{k_1,k_2,\cdots,k_{r-1}}(1-e^{1-e^t})\\
&=\sum_{l=0}^\infty \frac{t^l}{l!}\sum_{m=r-1}^\infty{m \brace k_1, k_2,\cdots, k_{r-1}}\frac{t^m}{m!}\\
&=\sum_{n=r-1}^\infty\bigg(\sum_{m=r-1}^n\binom{n}{m}{m \brace k_1,k_2,\cdots,k_{r-1}}\bigg)\frac{t^n}{n!}.
\end{split}
\end{equation}

On the other hand, by \eqref{eq18}, we get
\begin{equation}\label{eq25}
\begin{split}
\frac{d}{dt}\mathrm{Li}_{k_1,\cdots,k_{r-1},1}(1-e^{1-e^t})&=\frac{d}{dt}\sum_{n=r}^\infty{n \brace k_1,\cdots,k_{r-1},1}\frac{t^n}{n!}\\
&=\sum_{n=r}^\infty{n \brace k_1,\cdots,k_{r-1},1}\frac{t^{n-1}}{(n-1)!}\\
&=\sum_{n=r-1}^\infty{n+1 \brace k_1,\cdots,k_{r-1},1}\frac{t^n}{n!}.
\end{split}
\end{equation}

Therefore, by \eqref{eq24} and \eqref{eq25}, we obtain the following theorem.
\begin{theorem}

For $n\geq0$ and $r\geq0$ with $n+1\geq r$, we have
\begin{equation*}
\begin{split}
{n+1 \brace k_1,\cdots,k_{r-1},1}=\sum_{m=r-1}^n\binom{n}{m}{m \brace k_1,k_2,\cdots,k_{r-1}}.
\end{split}
\end{equation*}

\end{theorem}

Replacing $t$ by $e^t-1$ in \eqref{eq07}, we get
\begin{equation}\label{eq26}
\begin{split}
\frac{1}{(1-e^{1-e^t})^r}\mathrm{Li}_{k_1,\cdots,k_r}(1-e^{1-e^t})&=\sum_{m=0}^\infty B_m^{(k_1,\cdots,k_r)}\frac{1}{m!}(e^t-1)^m\\
&=\sum_{m=0}^\infty B_{m}^{(k_1,k_2,\cdots,k_r)}\sum_{n=m}^\infty {n \brace m}\frac{t^n}{n!}\\
&=\sum_{n=0}^\infty \bigg(\sum_{m=0}^nB_{m}^{(k_1,k_2,\cdots,k_r)}{n \brace m}\bigg)\frac{t^n}{n!}.
\end{split}
\end{equation}

On the other hand, by \eqref{eq18}, we get
\begin{equation}\label{eq27}
\begin{split}
&\frac{1}{(1-e^{1-e^t})^r}\mathrm{Li}_{k_1,k_2,\cdots,k_r}(1-e^{1-e^t})\\
&=\sum_{l=0}^\infty\binom{r+l-1}{l}e^{l(1-e^t)}\sum_{m=r}^\infty {m \brace k_1,\cdots, k_r}\frac{t^m}{m!}\\
&=\sum_{j=0}^\infty\bigg(\sum_{l=0}^\infty\binom{r+l-1}{l}\phi_j(-l)\bigg)\frac{t^j}{j!}\sum_{m=r}^\infty {m \brace k_1,\cdots, k_r}\frac{t^m}{m!}\\
&=\sum_{n=r}^\infty\bigg(\sum_{j=0}^{n-r}\sum_{l=0}^\infty \binom{r+l-1}{l}\phi_j(-l)\binom{n}{j}{n-j \brace k_1,\cdots, k_r}\bigg)\frac{t^n}{n!},
\end{split}
\end{equation}

where $\phi_n(x)$ are Bell polynomials given by
\begin{equation*}
\begin{split}
e^{x(e^t-1)}=\sum_{n=0}^\infty \phi_n(x)\frac{t^n}{n!}.
\end{split}
\end{equation*}

Therefore, by \eqref{eq26} and \eqref{eq27}, we obtain the following theorem.

\begin{theorem}

For $n,r\geq0$ with $n\geq r$, we have
\begin{equation*}
\begin{split}
\sum_{m=0}^nB_m^{(k_1,k_2,\cdots,k_r)}{n \brace m}=\sum_{l=0}^\infty\sum_{j=0}^{n-r}\binom{r+l-1}{l}\phi_j(-l)\binom{n}{j}{n-j \brace k_1,k_2,\cdots,k_r}.
\end{split}
\end{equation*}

In addition, for $0 \le n<r$, we have

\begin{equation*}
\begin{split}
\sum_{m=0}^n B_m^{(k_1,k_2,\cdots, k_r)}{n \brace m}=0.
\end{split}
\end{equation*}

\end{theorem}

In view of \eqref{eq11} and \eqref{eq14}, we may consider the unsigned multi-Stirling numbers of the first kind defined by
\begin{equation}\label{eq28}
\begin{split}
\mathrm{Li}_{k_1,k_2,\cdots,k_r}(t)=\sum_{n=r}^\infty{n \brack k_1,k_2,\cdots,k_r}\frac{t^n}{n!}.
\end{split}
\end{equation}

Note that
\begin{equation}\label{eq29}
\begin{split}
{n \brack \underbrace{1,1,\cdots,1}_{r-\text{times}}}={n \brack r},\quad (n,r\geq0).
\end{split}
\end{equation}

By replacing $t$ by $1-e^{1-e^t}$ in \eqref{eq28}, we set
\begin{equation}\label{eq30}
\begin{split}
\mathrm{Li}_{k_1,k_2,\cdots,k_r} (1-e^{1-e^t})&=\sum_{m=r}^\infty{{m \brack k_1,k_2,\cdots,k_r}}\frac{1}{m!}(1-e^{1-e^t})^m \\
&=\sum_{m=r}^\infty{m \brack k_1,k_2,\cdots,k_r}(-1)^m\sum_{l=m}^\infty {l \brace m}\frac{1}{l!}(1-e^t)^l \\
&=\sum_{l=r}^\infty \bigg(\sum_{m=r}^l (-1)^{l-m}{m \brack k_1,k_2,\cdots,k_r}{l \brace m}\bigg)\frac{1}{l!}(e^t-1)^l \\
&=\sum_{l=r}^\infty \bigg(\sum_{m=r}^l(-1)^{l-m} {m \brack k_1,k_2,\cdots,k_r}{l \brace m}\bigg)\sum_{n=l}^\infty {n \brace l}\frac{t^n}{n!} \\
&=\sum_{n=r}^\infty \bigg(\sum_{l=r}^n \sum_{m=r}^l (-1)^{l-m} {n \brace l}{l \brace m} {m \brack k_1,k_2,\cdots,k_r}\bigg)\frac{t^n}{n!}.
\end{split}
\end{equation}

Therefore, by \eqref{eq18} and \eqref{eq30}, we obtain the following theorem.

\begin{theorem}

For $n,r\in \mathbb{N}$ with $n\geq r$, we have
\begin{equation*}
\begin{split}
{n \brace k_1,k_2,\cdots,k_r} =\sum_{l=r}^n \sum_{m=r}^l (-1)^{l-m} {n \brace l}{l \brace m}{m \brack k_1,k_2,\cdots,k_r}.
\end{split}
\end{equation*}

\end{theorem}

Replacing $t$ by $e^t-1$ in \eqref{eq08}. Then we have
\begin{equation}\label{eq31}
\begin{split}
\frac{\mathrm{Li}_{k_1,k_2,\cdots,k_r} (1-e^{1-e^t})}{(2-e^{t})^r}&=\sum_{k=r}^\infty L^{(k_1,k_2,\cdots,k_r)}(k,r) \frac{1}{k!}(e^t-1)^k \\
&=\sum_{k=r}^\infty L^{(k_1,k_2,\cdots,k_r)}(k,r) \sum_{n=k}^\infty {n \brace k} \frac{t^n}{n!} \\
&=\sum_{n=r}^\infty \bigg(\sum_{k=r}^n {n \brace k}L^{(k_1,k_2,\cdots,k_r)}(k,r)\bigg)\frac{t^n}{n!}.
\end{split}
\end{equation}

On the other hand, by \eqref{eq15}, we get
\begin{equation}\label{eq32}
\begin{split}
\frac{1}{(2-e^t)^r}\mathrm{Li}_{k_1,k_2,\cdots,k_r} (1-e^{1-e^t})&=\sum_{k=r}^\infty {k \brace k_1,k_2,\cdots,k_r}\frac{t^k}{k!} \sum_{l=0}^\infty F_l^{(r)}(1) \frac{t^l}{l!} \\
&=\sum_{n=r}^\infty \bigg(\sum_{k=r}^n{k \brace k_1,k_2,\cdots,k_r}\binom{n}{k}F_{n-k}^{(r)}(1)\bigg)\frac{t^n}{n!}.
\end{split}
\end{equation}

Therefore, by \eqref{eq31} and \eqref{eq32}, we obtain the following theorem.
\begin{theorem}
For $n,r \in \mathbb{N}$ with $n\geq r$, we have
\begin{equation*}
\begin{split}
\sum_{k=r}^n{n \brace k}L^{(k_1,k_2,\cdots,k_r)}(k,r) =\sum_{k=r}^n {k \brace k_1,k_2,\cdots,k_r} \binom{n}{k} F_{n-k}^{(r)}(1).
\end{split}
\end{equation*}
\end{theorem}

Replacing $t$ by $\log(1+t)$ in \eqref{eq18}, we get
\begin{equation}\label{eq33}
\begin{split}
\mathrm{Li}_{k_1,k_2,\cdots,k_r} (1-e^{-t})&=\sum_{m=r}^\infty {m \brace k_1,k_2,\cdots,k_r} \frac{1}{m!} (\log(1+t))^m \\
&=\sum_{m=r}^\infty {m \brace k_1,k_2,\cdots,k_r} \sum_{n=m}^\infty S_1(n,m)\frac{t^n}{n!} \\
&=\sum_{n=r}^\infty \bigg(\sum_{m=r}^n {m \brace k_1,k_2,\cdots,k_r} S_1(n,m)\bigg)\frac{t^n}{n!}.
\end{split}
\end{equation}

On the other hand, by \eqref{eq08}, we get
\begin{equation}\label{34}
\begin{split}
\mathrm{Li}_{k_1,k_2,\cdots,k_r} (1-e^{-t})&=\frac{\mathrm{Li}_{k_1,k_2,\cdots,k_r} (1-e^{-t})}{(1-t)^r}(1-t)^r \\
&=\sum_{m=r}^\infty l^{(k_1,k_2,\cdots,k_r)}(m,r)\frac{t^m}{m!}\sum_{l=0}^\infty \binom{r}{l}(-1)^lt^l \\
&=\sum_{n=r}^\infty \bigg(\sum_{m=r}^n \frac{1}{m!}l^{(k_1,k_2,\cdots,k_r)}(m,r)\binom{r}{n-m}(-1)^{n-m}\bigg)t^n.
\end{split}
\end{equation}

Therefore, by \eqref{eq33} and \eqref{34}, we obtain the following theorem.
\begin{theorem}

For $n,r \in \mathbb{N}$ with $n\geq r$, we have
\begin{equation}
\begin{split}
\sum_{m=r}^n {n \brace k_1,k_2,\cdots,k_r} S_1(n,m) = \sum_{m=r}^n \frac{n!}{m!}l^{(k_1,k_2,\cdots,k_r)}(m,r) \binom{r}{n-m}(-1)^{n-m}.
\end{split}
\end{equation}

\end{theorem}

Replacing $t$ by $\log(1+t)$ in \eqref{eq18}, we get
\begin{equation}\label{eq35}
\begin{split}
\mathrm{Li}_{k_1,k_2,\cdots,k_r} (1-e^{-t})&=\sum_{m=r}^\infty {m \brace k_1,k_2,\cdots,k_r} \frac{1}{m!}(\log(1+t))^m \\
&=\sum_{n=r}^\infty{m \brace k_1,k_2,\cdots,k_r} \sum_{n=m}^\infty S_1(n,m)\frac{t^n}{n!} \\
&=\sum_{n=r}^\infty \bigg(\sum_{m=r}^n {m \brace k_1,k_2,\cdots,k_r}S_1(n,m)\bigg)\frac{t^n}{n!}.
\end{split}
\end{equation}

On the other hand, by \eqref{eq07}, we get
\begin{equation}\label{eq36}
\begin{split}
\mathrm{Li}_{k_1,k_2,\cdots,k_r} (1-e^{-t})&=\frac{\mathrm{Li}_{k_1,k_2,\cdots,k_r} (1-e^{-t})}{(1-e^{-t})^r}(1-e^{-t})^r \\
&=\bigg(\sum_{m=0}^\infty B_m^{(k_1,k_2,\cdots,k_r)}\frac{t^m}{m!}\bigg)r!(-1)^r\frac{1}{r!}(e^{-t}-1)^r \\
&=\bigg(\sum_{m=0}^\infty B_m^{(k_1,k_2,\cdots,k_r)}\frac{t^m}{m!}\bigg)r!(-1)^r \sum_{j=r}^\infty {j \brace r}(-1)^j \frac{t^j}{j!} \\
&=\sum_{n=r}^\infty \bigg(\sum_{j=r}^n {j \brace r}(-1)^{j-r}r! \binom{n}{j}B_{n-j}^{(k_1,k_2,\cdots,k_r)}\bigg) \frac{t^n}{n!}.
\end{split}
\end{equation}

Therefore, by \eqref{eq35} and \eqref{eq36}, we obtain the following theorem.

\begin{theorem}

For $n,r \in \mathbb{N}$ with $n\geq r$, we have
\begin{equation*}
\begin{split}
\frac{1}{r!}\sum_{m=r}^n {n \brace k_1,k_2,\cdots,k_r} S_1(n,m) = \sum_{j=r}^n {j \brace r}(-1)^{j-r} \binom{n}{j}B_{n-j}^{(k_1,k_2,\cdots,k_r)}.
\end{split}
\end{equation*}

\end{theorem}

\section{Conclusion}
With the help of the multiple logarithm we introduced the multi-Stirling numbers of the second kind which generalize the Stirling numbers of the second kind. Similarly, the multi-Stirling numbers of the first kind, the multi-Lah numbers and the multi-Bernoulli numbers, which respectively generalize the unsigned Stirling numbers of the first kind, the unsigned Lah numbers and the higher-order Bernoulli numbers, are defined by means of the multiple logarithm. In this paper, we derived several identities involving those four numbers defined by means of multiple logarithm and some other numbers. \par
The unsigned Lah numbers found their recent applications to the steganography method in telecommunication \cite{A} and the perturbative description of chromatic dispersion in optics \cite{B}. So it would be an interesting research problem to find some applications of the multi-Lah numbers. \par
As it has been, we would like to continue to find some special numbers and their applications to physics, science and engineering as well as to mathematics.

\bigskip

\noindent{\bf{Acknowledgments}} \\
All authors thank Jangjeon Institute for Mathematical Science for the support of this research.

\vspace{0.1in}

\noindent{\bf {Availability of data and material}} \\
Not applicable.

\vspace{0.1in}

\noindent{\bf{Funding}} \\
The first author of this research was supported by the Research Grant of Kwangwoon University
in 2023 and the third author was supported by the Basic Science Research Program, the National
               Research Foundation of Korea, (NRF-2021R1F1A1050151).

\vspace{0.1in}

\noindent{\bf{Ethics approval and consent to participate}} \\
All authors declare that there is no ethical problem in the production of this paper.

\vspace{0.1in}

\noindent{\bf {Competing interests}} \\
All authors declare no conflict of interest.

\vspace{0.1in}

\noindent{\bf{Consent for publication}} \\
All authors want to publish this paper in this journal.

\vspace{0.1in}

%


\end{document}